\documentclass[11pt]{article}
\usepackage{amssymb,latexsym}
\usepackage{euscript}
 \usepackage{epstopdf}
\usepackage{graphicx}
 \usepackage[applemac]{inputenc}
\usepackage{theorem}
\usepackage{epsfig}
\DeclareSymbolFont{AMSb}{U}{msb}{m}{n}
\DeclareSymbolFontAlphabet{\Bbb}{AMSb}

\newtheorem{theorem}{Theorem}[section]

{\theorembodyfont{\rmfamily}\newtheorem{defin}{Definition}[section]}
{\theorembodyfont{\rmfamily}\newtheorem{remq}{Remark}[section]}

\newcounter{example}[section]

\font\ly=lasy10 scaled 1100
\catcode`\«=\active\def«{\ofd}
\catcode`\»=\active\def»{\ffg}
\chardef\lg='050
\chardef\rg='051
\def\ofd{\leavevmode\hbox{\ly\lg\kern-0.2em\lg}}
\def\ffg{\leavevmode\hbox{\ly\rg\kern-0.2em\rg}}

\begin{document}
\title{Christoffel words and Markoff triples}

\author{
{\Large Christophe Reutenauer}
\thanks{Christophe Reutenauer; Universit\'e du Qu\'ebec \`a
Montr\'eal; Département de mathématiques; Case postale 8888,
succursale Centre-Ville, Montr\'eal (Qu\'ebec) Canada, H3C 3P8
(mailing address); e-mail: christo@math.uqam.ca.  Supported by NSERC.}\\ UQAM}
\date{\empty}
\maketitle

\setcounter{section}{0}
\section{Introduction}
A Markoff triple is a triple of natural integers $a,b,c$ which satisfies the Diophantine Equation $$a^2+b^2+c^2=3abc.$$ These numbers have been introduced in the work of Markoff \cite{M2}, where he finished his earlier work \cite{M1} on minima of quadratic forms and approximation of real numbers by continued fractions. He studies for this certain bi-infinite sequences and shows that they must be periodic; each of these sequences is a repetition of the same pattern, which turns out to be one of the words introduced some years before by Christoffel \cite{C1}; Markoff was apparently not aware of this work of Christoffel. These words were called Christoffel words much later in \cite{B} and may be constructed geometrically, see \cite{MH}, \cite{OZ}, \cite{B}, \cite{BL},  and have also many different interpretations, in particular in the free group on two generators (see \cite{BLRS} for the theory of Christoffel words).

In this Note, we describe a mapping which associates to each Christoffel word a Markoff triple and show that this mapping is a bijection. The construction and the result are a variant a theorem of Harvey Cohn \cite {C}, see also \cite{CF} chapter 7 and \cite{P}.

As them, we use the {\it Fricke identities}. We use also the tree construction of the Christoffel words, see \cite{BL} and \cite {BdL}, and moreover the arithmetic recursive construction of the solutions of the previous Diophantine equation, see \cite{M2}, \cite{F}, \cite{CF}, \cite{W}, \cite{P}.

%
%
%
%
\section{Results}
We consider {\it lattice paths}, which are consecutive elementary steps in the $x,y$-plane; each {\it elementary step} is a segment $[(a,b),(a+1,b)]$ or $[(a,b),(a,b+1)]$, with $a,b\in {\mathbf Z}$.

Let $p,q$ be relatively prime integers. Consider the segment from $(0,0)$ to $(p,q)$ and the lattice path from $(0,0)$ to $(p,q)$ located below this segment and such that the polygon delimited by the segment and the paths has no interior integer point.

The {\it Christoffel word of slope q/p} is the word in the free monoid $\{x,y\}^*$ coding the above path, where $x$ (resp. $y$) codes an horizontal (resp. vertical) elementary  step. See the figure, where is represented the path with $(p,q)=(7,3)$ corresponding to the Christoffel word of slope 3/7.

%
%
%
%
\begin{center}
\begin{figure}
[!ht]
{\centering
\includegraphics[width=2.3in]{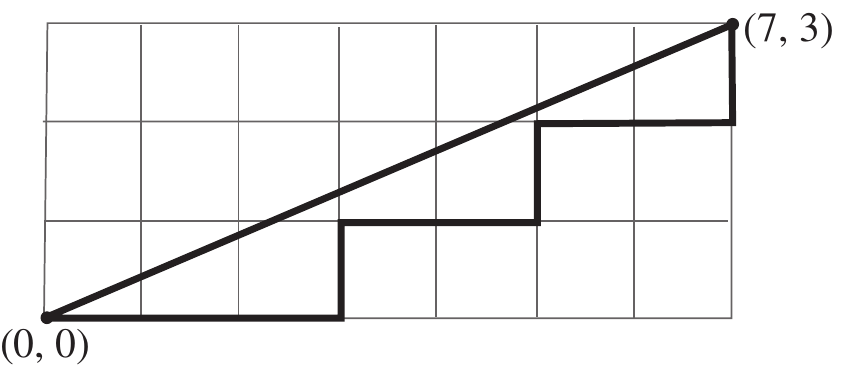}
\caption{the Christoffel word $xxxyxxyxxy$ of slope 3/7}
}
\end{figure}
\end{center}

Note that the definition includes the particular cases $(p,q)=(1,0)$ and $(0,1)$, corresponding to the Christoffel words $x$ and $y$. All other Christoffel words will be called {\it proper}. 

Each proper Christoffel word $w$ has a unique {\it standard factorization} $w=w_1w_2$; it is obtained by cutting the path corresponding to $w$ at the integer point closest to the segment. In the figure, the standard factorization is given by $w_1=xxxyxxy$, $w_2=xxy$. The words $w_1$ and $w_2$ are then Christoffel words.

Define the homomorphism $\mu$ from the free monoid $\{x,y\}^*$ into $SL_2(\mathbf{Z})$, defined by $\mu x=
\left( 
\begin{array}{cc}
2&1\\1&1
\end{array}
\right)
$
and $\mu y=
\left( 
\begin{array}{cc}
5&2\\2&1
\end{array}
\right)
$.

A {\it Markoff triple} is a multiset $\{a,b,c\}$ of positive integers which satisfies the equation $a^2+b^2+c^2=3abc$. The triple is {\it proper} if $a,b,c$ are distinct. It is classical that the only inproper Markoff triples are $\{1,1,1\}$ and $\{1,1,2\}$, see \cite{CF} Chapter 2. 

\begin{theorem}
A multiset $\{a,b,c\}$ is a proper Markoff triple if and only if it is equal to $\{1/3 Tr(\mu w_1), 1/3 Tr(\mu w_2), 1/3 Tr(\mu w)\}$ for some unique Christoffel word $w$ with standard factorization $w=w_1w_2$.
\end{theorem}

It is known that for each Christoffel word $w$, one has $1/3 Tr(\mu w) = \mu (w)_{1,2}$, see \cite{R} Lemme 3.2
or \cite{BLRS} Lemma 8.7.

{\bf Proof}
1. Each couple $(w_1,w_2)$ forming the standard factorization of a Christoffel word is obtained by applying iteratively the rules $(u,v) \to (u,uv)$ or $(u,v) \to (uv,v)$ starting from the couple $(x,y)$, see \cite {BL} Proposition 2 or \cite{BdL} Lemma 7.3. As a consequence, one obtains inductively that $w_1w_2w_1^{-1}w_2^{-1}$ is conjugate to $xyx^{-1}y^{-1}$ in the free group $F_2$ generated by $x$ and $y$.

Now, let  $w$ be some Christoffel word with standard factorization $w=w_1w_2$. Let $a=1/3Tr(w_1), b=1/3Tr(w_2), c=1/3Tr(w)$. We use the Fricke identity
$$ Tr(A)^2+Tr(B)^2+Tr(AB)^2=Tr(ABA^{-1}B^{-1})+2+Tr(A)Tr(B)Tr(AB) $$
for any $A,B\in SL_2(\mathbf{Z})$. We take $A=\mu (w_1)$, $B=\mu (w_2)$ and therefore $AB=\mu (w)$. 
Thus $9a^2+9b^2+9c^2=27abc$ (and we are done), provided $Tr(ABA^{-1}B^{-1})=-2$. Since $w_1w_2w_1^{-1}w_2^{-1}$ is conjugate to $xyx^{-1}y^{-1}$, it suffices to show that $Tr(\mu x\mu y\mu x^{-1}\mu y^{-1}) = -2$. Now the matrix $\mu x\mu y\mu x^{-1}\mu y^{-1}$ is equal to $
\left( 
\begin{array}{cc}
11&-24\\-6&-13
\end{array}
\right)
$, which shows what we want.

It remains to show that $a,b,c$ are distinct. In the special case where $w=xy$, this is seen by inspection. Furthermore, $\mu (w)_{12} = \mu (w_1)_{11}\mu (w_2)_{12}+\mu (w_1)_{12}\mu (w_2)_{22}$, which implies that $\mu (w)_{12}>\mu (w_1)_{12},\mu (w_2)_{12}$, since the matrices have positive coefficients. Thus, by the remark before the proof, $c>a,b$.  Now, since we assume $w\neq xy$ and by the tree construction recalled at the beginning of the proof, $w_1$ is a prefix of $w_2$ or $w_2$ is a suffix of $w_1$. Then $a<b$ or $a>b$ repectively.

2. Let  $\{a,b,c\}$ be a proper Markoff triple and assume that $a<b<c$. Then $\{a,b,3ab-c\}$ is a Markoff triple and $3ab-c<b$ by \cite{M2}, \cite{F}, \cite{CF}, \cite{W}.

Suppose first that this triple is proper. By induction, the multiset $\{a,b,3ab-c\}$ is equal to $\{1/3 Tr(\mu w_1), 1/3 Tr(\mu w_2), 1/3 Tr(\mu w)\}$ for some Christoffel word $w$ with standard factorization $w=w_1w_2$. By the first part of the proof, $b=1/3 Tr(\mu w)$, since both numbers are the maximum of their multiset. Then either (i) $a=1/3 Tr(\mu w_1)$ and $3ab-c=1/3 Tr(\mu w_2)$, or (ii) $a=1/3 Tr(\mu w_2)$ and $3ab-c=1/3 Tr(\mu w_1)$.

In case (i), we have $c=1/3 Tr(\mu (w_1^2 w_2))$; indeed, for $A,B$ in $SL_2(\mathbf Z)$, $Tr(A^2B)+Tr(B)=Tr(A)Tr(AB)$, hence $Tr(\mu (w_1^2 w_2))=Tr(\mu w_1)Tr(\mu w)-Tr(\mu w_2)=3a.3b-3(3ab-c)=3c$.

In case (ii), we have $c=Tr(\mu (w_1 w_2^2))$; indeed, $Tr(AB^2)+Tr(A)=Tr(AB)Tr(B)$, hence $Tr(\mu (w_1 w_2^2))=Tr(\mu w)Tr(\mu w_2)-Tr(\mu w_1)=3b.3a-3(3ab-c)=3c$.

Thus in case (i),  $\{a,b,c\}$ corresponds to the Christoffel word $w_1^2 w_2$ (with standard factorization $w_1.w_1w_2$) and in case (ii) to the Christoffel word $w_1 w_2^2$ (with standard factorization $w_1w_2.w_2$).

If $\{a,b,3ab-c\}$ is not proper, then, since $a\neq b$, we have $\{a,b,3ab-c\} = \{1,1,2\}$ and $a=1$ and $b=2$ or conversely. Thus $3ab-c=1$ and this implies $c=3.2.1-1=5$ and $\{a,b,c\}$ corresponds to the Christoffel word $xy$.

3. Concerning unicity, suppose that the proper Markoff triple $\{a,b,c\}$ with $a<b<c$ may be obtained from the two Christoffel words $u$ and $v$ with standard factorization $u=u_1u_2$ and $v=v_1v_2$. 
We may assume that they are both distinct from $xy$. Then by the tree construction of Christoffel words, $u_1$ is a prefix of $u_2$ or $u_2$ is a suffix of $u_1$, and similarly for $v_1$, $v_2$. Hence,  we are by symmetry reduced to  two cases: (i) $u_1$ is a prefix of $u_2$ and $v_1$ is a prefix of $v_2$, or (ii) $u_1$ is a prefix of $u_2$ and $v_2$ is a suffix of $v_1$.

In Case (i), $u_2$ (resp. $v_2$) has the standard factorization $u_2=u_1u'_2$ (resp.  $v_2=v_1v'_2$). Moreover $a=1/3Tr(\mu u_1), b=1/3Tr(\mu u_2), c=1/3Tr(\mu u)$ as we saw in the first part. Then the Markoff triple 
corresponding to the Christoffel word $u_1u'_2$ is equal to $\{a,3ab-c,b\}$. Indeed, $Tr(\mu u'_2)$
$ = Tr(\mu u_1)Tr(\mu (u_1u'_2))-Tr(\mu (u_1^2u'_2))=3a.3b-3c$, because $u_1u'_2=u_2$ and $u_1^2u'_2=u$. Similarly, the Markoff triple 
corresponding to the Christoffel word $v_1v'_2$
is also equal to $\{a,3ab-c,b\}$. Thus, by induction, the Christoffel words $u_1u'_2$ and $v_1v'_2$ are equal, together with their standard factorization, hence $u=v$.

In Case (ii), by the same calculation, the Markoff triple 
corresponding to the Christoffel word $u_1u'_2$
is equal to $\{a,3ab-c,b\}$, where the standard factorization of $u_2$ is $u_1u'_2$. Symmetrically, $v_1$ has the standard factorization $v_1=v'_1v_2$ and the Markoff triple 
corresponding to the Christoffel word $v'_1v_2$
is also equal to $\{a,3ab-c,b\}$. Thus, by induction, the Christoffel words $u_1u'_2$ and $v'_1v_2$ are equal, together with their standard factorization. Hence $u_1=v'_1$ and $u'_2=v_2$.  Thus $u=u_1u_2=u_1^2u'_2$ and $v=v_1v_2=v'_1v_2^2=u_1{u'_2}^2$. Thus we are reduced to the following particular case, which shows that Case (ii) cannot happen.

4. Let $w=uv$ be the standard factorization of the proper Christoffel word $w$. Then the Markoff triple corresponding to the Christoffel words $u^2v$ and $uv^2$ are distinct.

Indeed, disregarding the trivial case $w=xy$, we may assume that $u$ is a prefix of  $v$ or $v$ is a suffix of $u$. In the first case, each entry of the matrix $\mu u$ is stricly smaller than the corresponding entry of the matric $\mu v$. The same therefore holds for the matrices $\mu u\mu u\mu v$ and $\mu u\mu v\mu v$. Thus  $Tr(\mu (u^2v))$  is smaller than $Tr(\mu (uv^2))$, which proves that the  Markoff triples are distinct, since so are their greatest elements. The second case is similar.

\hfill{\raisebox{+.2mm}{\rule{.21cm}{.21cm}}

It is known that each Markoff number is of the form $1/3Tr(\mu w)$ for some Christoffel word $w$, see \cite{R} Cor.3.1 or \cite{BLRS} Th.8.4. It is however not known if the mapping associating to each Christoffel $w$ word the Markoff number $1/3Tr(\mu w)$ is injective. This is equivalent to the {\it Markoff numbers injectivity conjecture}, see \cite{F} page 614, \cite{CF}, \cite{W}.

%


\begin{thebibliography}{4444}



\bibitem[B]{B}
Berstel, J., Trac\'es de droites, fractions continues et distances discr\`etes, in: Mots, M\'elanges offerts \`a M.-P. Sch\"utzenberger, Herm\`es, Paris, 1990. 

\bibitem[BL]{BL}
Borel, J.-P., Laubie, F., {\it Quelques mots sur la droite projective réelle}, Journal de Théorie des Nombres de Bordeaux 5, 1993, 23-51.

\bibitem[BLRS]{BLRS}
Berstel, J., Lauve, A., Reutenauer, C., Saliola, F., Combinatorics on Words: Christoffel words and repetitions in words, CRM-AMS, 2008.

\bibitem[BdL]{BdL}
Berstel, J. and de Luca, A., Sturmian words, Lyndon words and trees, Theoretical Computer Science 178, 1997, 171-203.

\bibitem[C]{C}
Cohn, H., Markoff forms and primitives words, Mathematische Annalen 196, 1972, 8-22.

\bibitem[CF]{CF}
Cusick, T.W. and Flahive, M.E., The Markoff and Lagrange spectra,
American Mathematical Society, 1989.

\bibitem[C1]{C1}
Christoffel, E.B., {\it Observatio arithmetica}, Annali di Matematica 6, 1875, 145-152.



\bibitem[F]{F}
Frobenius, G.F., {\"U}ber die Markoffsche Zahlen, Sitzungsberichte der K{\"o}niglich Preussischen Akademie der Wissenschaften zu Berlin 1913, 458-487.

\bibitem[M1]{M1}
Markoff, A.A., {\it Sur les formes quadratiques binaires ind\'efinies},
Mathe\-matische Annalen 15, 1879, 381-496.

\bibitem[M2]{M2}
Markoff, A.A., {\it Sur les formes quadratiques binaires ind\'efinies
(second m\'emoire)}, Mathematische Annalen 17, 1880, 379-399.


\bibitem[MH]{MH}
Morse, M. and Hedlund, G.A., {\it Symbolic dynamics II: Sturmian
Trajectories}, Amer. J. Math. 62, 1940, 1-42. 

\bibitem[OZ]{OZ}
Osborne, R.P., Zieschang, H., Primitives in the free group on two generators, Inventiones Mathematicae 63, 1981, 17-24.

\bibitem[P]{P}
Perrine, S., L'interpr\'etation matricielle de la th\'eorie de Markoff classique, International Journal of mathematical Sciences, 32, 2002, 193-262.

\bibitem[R]{R}
Reutenauer, C., Mots de Lyndon g\'en\'eralis\'es, S\'eminaire Lotharingien de Combinatoire 54, 16 pages, 2006.


\bibitem[W]{W}
Waldschmidt, M., Open Diophantine problems, Moscow mathematical Journal 4, 2004, 245-305.

\end{thebibliography}
\end{document}